\documentclass[10pt]{amsart}   	
									
\usepackage{amsmath}
\usepackage{amssymb, amsfonts, amsrefs}
\usepackage{amsthm}
\usepackage{hyperref}
\usepackage{mathtools}
\usepackage{tikz}
\usetikzlibrary{cd}
\usetikzlibrary{decorations.markings,decorations.pathreplacing,arrows}
\usepackage{enumitem}
\usepackage{makecell}
\setcellgapes{4pt}

\usepackage[margin=1in]{geometry}

\newtheorem{theorem}{Theorem}[section]
\newtheorem{lemma}[theorem]{Lemma}
\newtheorem{corollary}[theorem]{Corollary}
\newtheorem{proposition}[theorem]{Proposition}

\theoremstyle{definition}
\newtheorem{definition}[theorem]{Definition}
\newtheorem{example}[theorem]{Example}
\newtheorem{notation}[theorem]{Notation}

\renewcommand{\pmod}[1]{\textnormal{(mod #1)}}

\usepackage{soul}

\title{A Counterexample to the $\phi$-Dimension Conjecture}
\author{Eric J. Hanson}\address{Brandeis University, Department of Mathematics, 415 South Street, Waltham MA 02453, USA}\email{ehanson4@brandeis.edu}
\author{Kiyoshi Igusa}\address{Brandeis University, Department of Mathematics, 415 South Street, Waltham MA 02453, USA}\email{igusa@brandeis.edu}



\begin{document}

\noindent This version of the article has been accepted for publication, after peer review (when applicable) and is subject to Springer Nature's \href{https://www.springernature.com/gp/open-research/policies/accepted-manuscript-terms}{AM terms of use}, but is not the Version of Record and does not reflect post-acceptance improvements, or any corrections. The Version of Record is available online at: https://doi.org/10.1007/s00209-021-02795-7.

\bigskip

\maketitle

\begin{abstract}
	In 2005, the second author and Todorov introduced an upper bound on the finitistic dimension of an Artin algebra, now known as the $\phi$-dimension. The $\phi$-dimension conjecture states that this upper bound is always finite, a fact that would imply the finitistic dimension conjecture. In this paper, we present a counterexample to the $\phi$-dimension conjecture and explain where it comes from. We also discuss implications for further research and the finitistic dimension conjecture.
\end{abstract}

\section{Introduction}
	Let $\Lambda$ be an Artin algebra; that is, an associative algebra which is finitely generated as a module over a commutative Artinian ring. We denote by $\mathsf{mod}\Lambda$ the category of finitely generated (left) $\Lambda$-modules and by $\mathsf{proj}\Lambda$ the full subcategory of finitely generated projective modules. Let $\Omega_{\Lambda}:\mathsf{mod}\Lambda\rightarrow\mathsf{mod}\Lambda$ be the syzygy functor. That is, on objects, $\Omega_\Lambda M$ is the kernel of the projective cover $M\xleftarrow{q} P$.
	
	One of the fundamental homological invariants is the \emph{projective dimension} of a module. This can be defined as
	$$\mathsf{pd}_\Lambda M := \min\left(\left\{k \in \mathbb{N}\mid\Omega_\Lambda^k M\in\mathsf{proj}\Lambda\right\}\cup\{\infty\}\right).$$
	This leads to the definitions of the \emph{projective dimension} and \emph{finitistic dimension} of an algebra, given by
	\begin{eqnarray*}
		\mathsf{gl.dim}\Lambda &:=& \sup\left\{\mathsf{pd}_\Lambda M\mid M \in \mathsf{mod}\Lambda\right\} = \max\{\mathsf{pd}_\Lambda S\mid S\in\mathsf{mod}\Lambda \textnormal{ is simple}\},\\
		\mathsf{fin.dim}\Lambda &:=& \sup\left\{\mathsf{pd}_\Lambda M\mid M \in \mathsf{mod}\Lambda,\mathsf{pd}_{\Lambda}M < \infty\right\}.
	\end{eqnarray*}
	
	It is clear that if $\mathsf{gl.dim}\Lambda < \infty$ then $\mathsf{gl.dim}\Lambda = \mathsf{fin.dim}\Lambda$, but it remains unknown whether the converse is true. The (small) finitistic dimension conjecture, first formally communicated by Bass in 1960 \cite{bass_finitistic}, states precisely that this is the case. In other words, it states that $\mathsf{fin.dim}\Lambda < \infty$ for all Artin algebras.
	
	This conjecture has motivated an enormous quantity of research since it was first stated, both for its intrinsic interest and for its relation to other open questions in representation theory. For example, it is known that the finitistic dimension conjecture implies the generalized Nakayama conjecture (posed by Auslander-Reiten \cite{auslander_generalized}), the Wakamatsu tilting conjecture (posed by Beligiannis-Reiten \cite[Chapter 3]{beligiannis_homological}), and the Gorenstein symmetry conjecture (noted by Auslander). A more thorough description of the relationship between these conjectures can be found in \cite{kirkman_finitistic}, \cite{xi_finitistic} and \cite{wei_finitistic2}.
	
	To study the finitistic dimension conjecture, the second author and Todorov introduced two new homological invariants, which they denoted $\phi$ and $\psi$, in \cite{igusa_finitistic}. From these invariants, which are now sometimes referred to as the \emph{Igusa-Todorov functions}, one defines the \emph{$\phi$-dimension and $\psi$-dimension} of an Artin algebra as
	\begin{eqnarray*}
		\phi\mathsf{dim}\Lambda &:=& \sup\{\phi(M)\mid M \in \mathsf{mod}\Lambda\},\\
		\psi\mathsf{dim}\Lambda &:=& \sup\{\psi(M)\mid M \in \mathsf{mod}\Lambda\}.
	\end{eqnarray*}
	These invariants have remained an active area of research since their inception, mainly for their relationship to other homological measures and the finitistic dimension conjecture. For example, in \cite{elsener_syzygies} and \cite{lanzilotta_igusa}, Elsener-Schifler and Lanzilotta-Mata independently prove that for Gorenstein algebras the $\phi$-dimension, $\psi$-dimension, and Gorenstein dimension all coincide. In addition, in \cite{huard_self}, Huard-Lanzilotta show that an algebra has $\phi\mathsf{dim}\Lambda = 0$ if and only if it is self-injective. It remains an open question to characterize algebras with $\phi\mathsf{dim}\Lambda = 1$. \c{S}en shows in \cite{sen_phi} that the $\phi$-dimension of any cyclic Nakayama algebra of infinite global dimension is even.
	
	The classes of algebras that have been shown to satisfy the finitistic dimension conjecture using the $\phi$ and $\psi$ functions include: algebras of representation dimension at most 3 \cite{igusa_finitistic}, algebras of finite injective dimension \cite{lanzilotta_igusa}, Gorenstein algebras \cite{lanzilotta_igusa}, truncated path algebras \cite{barrios_igusa}, monomial relation algebras \cite{lanzilotta_igusa}, Igusa-Todorov algebras \cite{wei_finitistic}, and, indirectly, special biserial algebras \cite{erdmann_radical}.
	
	As in \cite{lanzilotta_igusa}, we recall that
	$$\mathsf{fin.dim}\Lambda \leq \phi\mathsf{dim}\Lambda \leq \psi\mathsf{dim}\Lambda \leq \mathsf{gl.dim}\Lambda.$$
	In particular, finiteness of either the $\phi$- or $\psi$-dimension implies finiteness of the finitistic dimension. This fact, and the prevalence of the $\phi$- and $\psi$-dimensions in recent literature on the finitistic dimension conjecture, has led to the so called $\phi$-dimension conjecture and $\psi$-dimension conjecture, formally stated by Fernandes-Lanzilotta-Mendoza \cite{fernandes_phi} and Lanzilotta-Mendoza \cite{lanzilotta_relative}. These conjectures state, respectively, that $\phi\mathsf{dim}\Lambda < \infty$ and $\psi\mathsf{dim}\Lambda < \infty$ for all Artin algebras. We observe that the $\phi$-dimension conjecture implies the finitistic dimension conjecture. Said another way, if there is a counterexample to the finitistic dimension conjecture, it must have infinite $\phi$-dimension. The purpose of this paper is to show that Artin algebras of infinite $\phi$-dimension exist, disproving the $\phi$-dimension conjecture (and hence also the $\psi$-dimension conjecture), while also showing that this class of algebras is important for future research and the search for a resolution of the finitistic dimension conjecture. Recently, Barrios-Mata have given another example of an Artin algebra of infinite $\phi$-dimension in \cite{barrios_algebras} which is independent of this paper.
	
\subsection{Our Counterexample}

We fix for the remainder of this paper an arbitrary field $K$. For example, we could take $K = \mathbb{C}$. Let $A = KQ/\mathsf{rad}^2 KQ$ and $A_3^{CT} = KC_3/\mathsf{rad}^2 KC_3$, where $Q$ and $C_3$ are the quivers shown below.

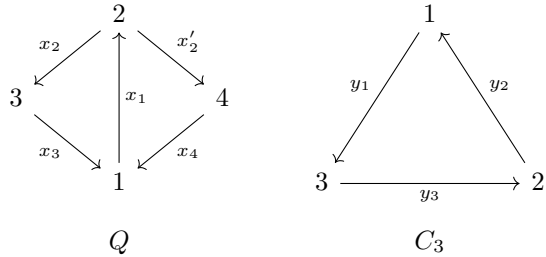
\begin{center}
\begin{tikzcd}[row sep = tiny]
		&2\arrow[dl,"x_2" above left]\arrow[dr,"x_2'"above right]&&&1\arrow[ddl,"y_1"above left]\\
		3\arrow[dr,"x_3"below left]&&4\arrow[dl,"x_4"below right]&\\
		&1\arrow[uu,"x_1"right]&&3\arrow[rr,"y_3"above]&&2\arrow[uul,"y_2"above right]\\[-.8em]
		&Q&&&C_3\\
\end{tikzcd}
\end{center}

That is, $A$ and $A_3^{CT}$ are the quotients of the path algebras (over $K$) of these quivers modulo the relations that the composition of any two arrows is zero. Readers unfamiliar with path algebras and their quotients are referred to \cite[Chapter II, III]{assem_elements}. 

For emphasis, we restate the following result of Lanzilotta-Mata.

\begin{theorem}\cite[Corollary 3.9]{lanzilotta_igusa}
	Let $\Lambda$ be a monomial relation algebra. Then $\phi\mathsf{dim} \Lambda < \infty$.
\end{theorem}

In particular, both $A$ and $A_3^{CT}$ are monomial relation algebras, so their $\phi$-dimensions are both finite. The counterexample we present in this paper is the algebra $A\otimes_K A_3^{CT}$. That is, the aim of this paper is to prove the following.

\begin{theorem}[Theorem \ref{thm:phidiminfinite}]
	The $\phi$-dimension of $A\otimes_K A_3^{CT}$ is infinite.
\end{theorem}

As the algebra $A\otimes_K A_3^{CT}$ is not a monomial relation algebra, this result is not in conflict with that of Lanzilotta-Mata. We also note that the algebra $A\otimes_K A_3^{CT}$ is not a counterexample to the finitistic dimension conjecture. Indeed, as a consequence of \cite[Thorem 16]{eilenberg_dimension}, we know $\mathsf{fin.dim}\left(A\otimes_K A_3^{CT}\right) = 0$. The authors thank Rene Marczinzik for pointing out this fact.

\subsection{Organization}
	
This paper is organized as follows: In Section \ref{sec:amalgamation}, we discuss the larger context from which we discovered this counterexample. In Section \ref{sec:phi}, we recall the definitions of the $\phi$ and $\psi$ functions and some useful lemmas. In Section \ref{sec:periodic}, we study the category of \emph{3-periodic chain complexes} of an Artin algebra over a field. This category is isomorphic to the module category of the algebra formed by taking a tensor product with $A_3^{CT}$, and our counterexample is a specific instance of this construction. In Section \ref{sec:syzygy}, we give a procedure for computing the syzygies of a class of objects, which we call \emph{truncated projective resolutions}, in the category of (bounded) chain complexes and the category of 3-periodic chain complexes. In Section \ref{sec:computation}, we give an explicit computation which shows that the $\phi$-dimension of our proposed counterexample is infinite (i.e., we prove Theorem \ref{thm:phidiminfinite}) using the results of Sections \ref{sec:periodic} and \ref{sec:syzygy}.

\section{A Brief Description of Amalgamation}\label{sec:amalgamation}

The goal of this section is explain the larger context in which we discovered this counterexample. The ideas explained here are a part of a larger forthcoming work by Gordana Todorov and the two authors of this paper. This was originally motivated by the joint work of the second author and Daniel \'Alvarez-Gavela, who use amalgamation to describe examples and invariants in contact topology in \cite{alvarez_turaev}.

In \cite{fock_cluster}, Fock and Goncharov define a combinatorial method for gluing together two cluster $\mathcal{X}$-varieties, which they refer to as \emph{amalgamation}. This process is described by how it behaves on seeds. In \cite{arkani_scattering}, amalgamation is used to study on-shell diagrams with the goal of understanding scattering amplitudes and the Jacobian algebra of a plabic diagram. As part of this study, an interpretation of amalgamation in terms of quivers containing ``half-arrows" is given.

In our larger work, we avoid the notion of half arrows by instead adding formal inverses to the arrows of the quivers and relations to make these new arrows redundant. The reason we do this is that after adding these `redundant' arrows, amalgamating, and deleting the remaining redundant arrows, we recover the Jacobian algebra of the amalgamation defined by Fock-Goncharov. This new interpretation also allows us to naturally expand the process of amalgamation to quivers with relations, and thus to quotients of path algebras. In our larger work, we study what happens when two algebras are amalgamated. We also define a pseudo-inverse process we call ``unamalgamation'' and study what happens when an algebra is unamalgamated.

The technical definition of amalgamation will not be a part of this paper in favor of an example. Consider the following pair of identical quivers of type $C_3$

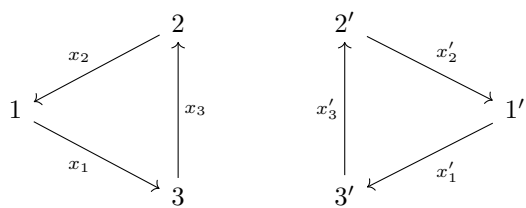
\begin{center}
	\begin{tikzcd}[row sep = tiny]
		&&2\arrow[dll,"x_2" above left]&&2'\arrow[drr,"x_2'"above right]\\
		1\arrow[drr,"x_1"below left]&&&&&&1'\arrow[dll,"x_1'"below right]\\
		&&3\arrow[uu,"x_3"right]&&3'\arrow[uu,"x_3'"left]&\\
	\end{tikzcd}
\end{center}
modulo the relations $\mathsf{rad}^2 C_3 = 0$ on each piece. Observe that this corresponds to two algebras of the form $A_3^{CT}$. We can then add formal inverses of each arrow, with relations to make them redundant. The result is the pair of quivers
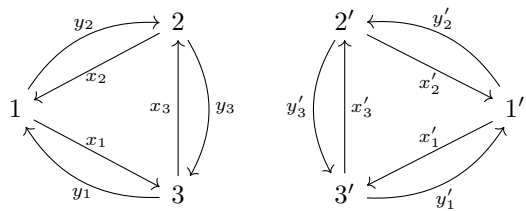
\begin{center}
	\begin{tikzcd}[row sep = tiny]
		&&2\arrow[dll,"x_2" below]\arrow[dll,leftarrow, bend right,"y_2" above]&&2'\arrow[drr,"x_2'"below]\arrow[drr,leftarrow,bend left,"y'_2"above]\\
		1\arrow[drr,"x_1"above]\arrow[drr,leftarrow, bend right,"y_1"below]&&&&&&1'\arrow[dll,"x_1'"above]\arrow[dll,leftarrow,bend left,"y'_1"below]\\
		&&3\arrow[uu,"x_3"left]\arrow[uu,leftarrow,bend right,"y_3"right]&&3'\arrow[uu,"x_3'"right]\arrow[uu,leftarrow,bend left,"y'_3"left]&\\
	\end{tikzcd}
\end{center}
modulo the relations from before and the new relations that $y_i = x_{i-1}x_{i-2}$ and likewise for $y'_i$ (where indices are considered mod 3). We observe that in this case, this means that each $y_i$ is 0, but if we had started with different relations this may not be the case.

The idea of amalgamation is to now identify an arrow of the first quiver with an arrow of the second. Up to symmetry, there are two ways we can do this.

Our first choice is to identify $x_3$ with $y_3' = x'_2x'_1$ and $x'_3$ with $y_3 = x_2x_1$ (hence identifying the vertex 3 with 2' and the vertex 2 with 3'). Upon making this identification and deleting unnecessary arrows, we are left with the quiver
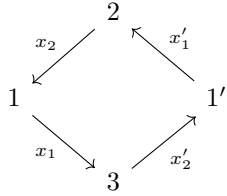
\begin{center}
	\begin{tikzcd}[row sep = tiny]
		&2\arrow[dl,"x_2" above left]\arrow[dr,leftarrow,"x_1'"above right]\\
		1\arrow[dr,"x_1"below left]&&1'\arrow[dl,leftarrow,"x_2'"below right]\\
		&3\\
	\end{tikzcd}
\end{center}
modulo the relations $x_2x_1 = 0$ and $x_2'x_1' = 0$.

Our other choice is to identify $x_3$ with $x'_3$ and $y_3$ with $y'_3$ (hence identifying the vertex 3 with 3' and the vertex 2 with 2'). Upon making this identification and deleting unnecessary arrows, we are left with the quiver

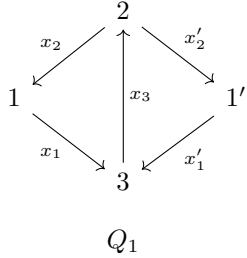
\begin{center}
	\begin{tikzcd}[row sep = tiny]
		&2\arrow[dl,"x_2" above left]\arrow[dr,"x_2'"above right]\\
		1\arrow[dr,"x_1"below left]&&1'\arrow[dl,"x_1'"below right]\\
		&3\arrow[uu,"x_3"right]\\
		&Q_1
	\end{tikzcd}
\end{center}
modulo the relation $\mathsf{rad}^2 Q_1 = 0$. We observe that upon taking a tensor product with $A_3^{CT}$, this quiver with relations is precisely our counterexample. Thus, in an informal sense, our counterexample is formed by putting together three copies of $A_3^{CT}$, a very well-understood algebra.

We now give a brief description of another problem we have studied using amalgamation (see \cite{hanson_resolution}).

\subsection{Global Dimension and Amalgamation}

Recall that a connected algebra $\Lambda = KQ/I$ is a \emph{Nakayama} algebra if and only if $Q$ is one of the quivers
	\begin{center}
		\begin{tikzpicture}
			\node at (0,0) {1};
			\node at (1,0) {2};
			\node at (3.5,0) {$n$};
			\draw [<-](.2,0)--(.8,0) node[midway,anchor=north]{$\gamma_1$};
			\draw [<-](1.2,0)--(1.8,0) node[midway,anchor=north]{$\gamma_2$};
			\node at (2.25,0) {$\cdots$};
			\draw [<-](2.7,0)--(3.3,0) node[midway,anchor=north]{$\gamma_{n-1}$};
			\draw [<-](3.5,.2)--(3.5,.5)--(0,.5)--(0,.2);
			\draw (0,.5)--(3.5,.5) node[midway,anchor=south]{$\gamma_n$};
			\node at (1.75,-1){$\Delta_n$};
		\begin{scope}[shift = {(5,0)}]
			\node at (0,0) {1};
			\node at (1,0) {2};
			\node at (3.5,0) {$n$};
			\draw [<-](.2,0)--(.8,0) node[midway,anchor=north]{$\gamma_1$};
			\draw [<-](1.2,0)--(1.8,0) node[midway,anchor=north]{$\gamma_2$};
			\node at (2.25,0) {$\cdots$};
			\draw [<-](2.7,0)--(3.3,0) node[midway,anchor=north]{$\gamma_{n-1}$};
			\node at (1.75,-1){$A_n$};
		\end{scope}
		\end{tikzpicture}
	\end{center}
	and $I$ is an admissible ideal. Nakayama algebras with quiver $\Delta_n$ are sometimes referred to as \emph{cyclic Nakayama algebras} or \emph{cycle algebras}.
	
	In the example above, we can consider amalgamation as a prodecure for combining two Nakayama algebras. Identifying two arrows in the opposite direction (the first case) resulted in another Nakayama algebra and identifying two arrows in the same direction (the second case) did not.
	
	Much work has been done to determine when a Nakayama algebra is of infinite global dimension. We now recall two approaches to this problem and the relationship between them.
	
	In \cite{ringel_gorenstein}, Ringel introduces the \emph{resolution quiver} of a Nakayama algebra as a tool to study its Gorenstein projective modules and determine when it is Gorenstein. Shen then shows in \cite{shen_homological} that a Nakayama algebra has finite global dimension if and only if its resolution quiver is connected and has \emph{weight} 1, an easy condition to check.
	
	In \cite{igusa_cyclic}, the second author and Zacharia study a larger class of algebras called \emph{monomial relation algebras}; that is, algebras of the form $KQ/I$ where $I$ is an admissible ideal generated by monomials and $Q$ is arbitrary. Their main result is that a monomial relation algebra has finite global dimension if and only if the \emph{cyclic homology} of its radical is trivial. They further show that it is sufficient to study only the cyclic homology of the \emph{overlying cycle algebras}, that is, cyclic Nakayama algebras $K\Delta_n/I'$ for which there is a map of quivers $\Delta_n \rightarrow Q$ and $I'$ is the pullback of $I$ by this map. For example, $KQ_1/\mathsf{rad}^2 KQ_1$ defined above in our example has $K\Delta_6/\mathsf{rad}^2K\Delta_6$ as an overlying Nakayama algebra coming from the map of quivers sending the cycle (1,2,3,4,5,6) to the path (3,2,1,3,2,1').
	
	In the case that a monomial relation algebra is Nakayama, the second author and Zacharia test whether the cyclic homology is trivial by relating it to the homology of a simplicial complex defined from the relations of the algebra, which we call the \emph{relation complex}. The main result of our concurrent work \cite{hanson_resolution} is to show that the Euler characteristic of the relation complex is exactly the number of connected components with weight 1 in the resolution quiver. Our proof of this uses (un)amalgamation.
	
	This result drove us to look for a generalization of resolution quivers to arbitrary monomial relation algebras. Made more precise, our results show that resolution quivers give information about the cyclic homology of Nakayama algebras, but the results of the second author and Zacharia apply to all monomial relation algebras. We wish to recover this result in the more general setting from a construction that restricts to a resolution quiver in the Nakayama case. As a starting point, we examined algebras formed by amalgamation. We discovered that the algebra $KQ'/I'$, where $Q'$ is the quiver
	
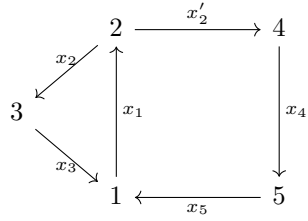
\begin{center}
	\begin{tikzcd}[row sep = tiny]
		&2\arrow[dl,"x_2" above]\arrow[rr,"x_2'"above]&&4\arrow[dd,"x_4"right]\\
		3\arrow[dr,"x_3"below]\\
		&1\arrow[uu,"x_1"right]&&5\arrow[ll,"x_5"below]\\
	\end{tikzcd}
\end{center}
contains overlying cycle algebras of arbitrarily high $\phi$-dimension. In an attempt to leverage this fact, we examined what happens to a (truncated) projective resolution for an overlying cycle algebra when we project it back to the original algebra. These resolutions can be considered as modules over a finite dimensional algebra by taking a tensor product with $A_3^{CT}$ (this idea is detailed in Section \ref{sec:periodic}). We then discovered that after contracting the arrow $x_4$, we could use this approach to generate modules of arbitrarily high $\phi$-dimension.

\section{The $\phi$ and $\psi$ Functions}\label{sec:phi}
	The aim of this section is to recall the definitions of the $\phi$ and $\psi$ functions, first given by the second author and Todorov, and prove a useful lemma. In this section, we let $\Lambda$ be an arbitrary Artin algebra; however, in the sections to follow we will restrict to the case where $\Lambda$ is an algebra over the field $K$.
	
\begin{definition}
	Fix a representative $[M]$ for every isomorphism class in $\mathsf{mod}\Lambda$. Let $K_0(\Lambda)$ be the quotient of the free abelian group generated by the symbols $[M]$ modulo the subgroup generated by
	\begin{itemize}
		\item $[M] - [N] - [N']$ for $M \cong N\oplus N'$
		\item $[P]$ for $P \in \mathsf{proj}\Lambda$.
	\end{itemize}
	$K_0(\Lambda)$ is called the \emph{split Grothendieck group} of $\Lambda$. It is well known that $K_0(\Lambda)$ is the free abelian group generated by the symbols $[M]$ for all isomorphism classes of indecomposable non-projective modules.
\end{definition}

We recall that the syzygy functor $\Omega_\Lambda$ induces a homomorphism $L:K_0(\Lambda)\rightarrow K_0(\Lambda)$ by defining $L[M] := [\Omega_\Lambda M]$. For $M \in \mathsf{mod}\Lambda$, let $\langle\mathsf{add} M\rangle$ be the subgroup of $K_0(\Lambda)$ generated by the direct summands of $M$. Observe that $\langle\mathsf{add} M\rangle$, and more generally $L^t\langle\mathsf{add} M\rangle$ for $t \in \mathbb{N},$ is also free abelian and thus has a well defined rank.

\begin{definition}\cite{igusa_finitistic}
	Let $M \in \mathsf{mod}\Lambda$. Then define
	\begin{eqnarray*}
		\phi(M) &:=& \min\{t\mid \forall j \in \mathbb{N}: \mathsf{rank}(L^t\langle\mathsf{add} M\rangle) = \mathsf{rank} (L^{t+j}\langle\mathsf{add} M\rangle)\}\\
		\psi(M) &:=& \phi(M) + \sup\{\mathsf{pd}_\Lambda N\mid N \textnormal{ a direct summand of } \Omega_\Lambda^{\phi(M)}M,\mathsf{pd}_\Lambda N < \infty\}.
	\end{eqnarray*}
\end{definition}
It is immediate that $\phi(M) < \infty$ for all $M \in \mathsf{mod}\Lambda$ and that if $\mathsf{pd}_\Lambda M < \infty$ then $\phi(M) = \psi(M) = \mathsf{pd}_\Lambda M$. We recall from the introduction that the $\phi$- and $\psi$-dimensions are then defined as
	\begin{eqnarray*}
		\phi\mathsf{dim}\Lambda &:=& \sup\{\phi(M)\mid M \in \mathsf{mod}\Lambda\}\\
		\psi\mathsf{dim}\Lambda &:=& \sup\{\psi(M)\mid M \in \mathsf{mod}\Lambda\}
	\end{eqnarray*}
We conclude this section with the following lemma, which will be used to prove the main theorem of this paper.

\begin{lemma} \label{lem:bound}Let $M, N \in \mathsf{mod}\Lambda$.
	\begin{enumerate}
		\item If $[M] = [N]$ in $K_0(\Lambda)$, then $\Omega_\Lambda M \cong \Omega_\Lambda N$.
		\item Suppose there exists $t \in \mathbb{N}$ so that $\Omega_\Lambda^t M \cong \Omega_\Lambda^t N$ and $\Omega_\Lambda^{t-1} M \ncong \Omega_\Lambda^{t-1} N$. Then $\phi\mathsf{dim}(M\oplus N) \geq t-1$.
	\end{enumerate}
\end{lemma}

\begin{proof}
	(1) Suppose $[M] = [N]$. Then there exists $R \in \mathsf{mod}\Lambda$ with no projective direct summands and $P,Q \in\mathsf{proj}\Lambda$ so that $M \cong R \oplus P$ and $N \cong R \oplus Q$. We then have that $\Omega_\Lambda M \cong \Omega_\Lambda R \cong \Omega_\Lambda N$.
	
	(2) Let $a = [M]-[N] \in K_0(\Lambda)$. By assumption, we have $a \in \langle\mathsf{add} (M \oplus N)\rangle$. Now since $\Omega_\Lambda^{t-1}M \ncong \Omega_\Lambda^{t-1}N$, we have by (1) that $[\Omega_\Lambda^{t-2}M] \neq [\Omega_\Lambda^{t-2}N]$; that is, $L^{t-2}(a) \neq 0$. Moreover, since $\Omega_\Lambda^t M \cong \Omega_\Lambda^t N$, we have $[\Omega_\Lambda^tM] = [\Omega_\Lambda^tN]$; that is, $L^t(a) = 0$. This means $\mathsf{rank}\left(L^{t}\langle\mathsf{add} M\oplus N\rangle\right) < \mathsf{rank}\left(L^{t-2}\langle\mathsf{add} M\oplus N\rangle\right)$, which implies the result.
\end{proof}


\section{The Category of 3-Periodic Chain Complexes}\label{sec:periodic}
	For the remainder of this paper, the term \emph{algebra} will be used to mean an elementary algebra over the field $K$. We recall that if $\Lambda$ is such an algebra, then we can express $\Lambda = KQ/I$ as the quotient of the ($K$-)path algebra of some quiver $Q$ by some admissible ideal $I$.

	Recall from the introduction that $A_3^{CT} = KC_3/\mathsf{rad}^2{KC_3}$, where $C_3$ is an oriented 3-cycle:
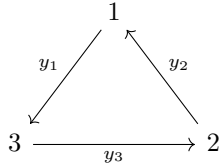
\begin{center}
\begin{tikzcd}[row sep = tiny]
		&1\arrow[ddl,"y_1"above left]\\
		\\
		3\arrow[rr,"y_3"below]&&2\arrow[uul,"y_2"above right]\\[-1em]
\end{tikzcd}
\end{center}
We remark that this algebra is cluster-tilted of type $A_3$, which is our reason for denoting it by $A_3^{CT}$. Now let $\Lambda = KQ/I$ be an arbritrary algebra. This section is devoted to the study of algebras of the form $\Lambda \otimes_K A_3^{CT}$.

It is well known that the quiver of $\Lambda\otimes_K A_3^{CT}$ is $Q\times C_3$ (see for example \cite{herschend_tensor}). For simplicity, we write this quiver in the form 
	
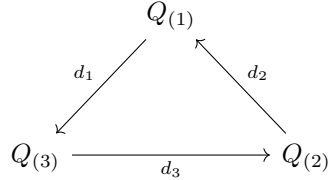
\begin{center}
	\begin{tikzcd}[row sep = tiny]
		&Q_{(1)}\arrow[ddl,"d_1"above left]\\
		\\
		Q_{(3)}\arrow[rr,"d_3"below]&&Q_{(2)}\arrow[uul,"d_2"above right]\\[-1em]
	\end{tikzcd}
	\end{center}
where by abuse of notation, each $d_i$ refers to $|Q_0|$ distinct arrows. More precisely, for all $v \in Q_0$ a vertex of $Q$, $d_1$ refers to the arrow $(v,1) \rightarrow (v,3)$ in $Q\times C_3$ and likewise for $d_2$ and $d_3$. There are then three types of relations:
	\begin{enumerate}
		\item For $R$ a relation of $A$, there is a relation $(R,1_{A_3^{CT}})$, where $1_{A_3^{CT}}$ is the identity of $A_3^{CT}$.
		\item Abusing notation, there is a relation $d^2 = 0$. This is the relation inherited from $A_3^{CT}$.
		\item For all arrows $\gamma: s \rightarrow t$ in $Q$, there is a relation $\gamma\circ d = d\circ \gamma$ (again abusing notation).
	\end{enumerate}
	
Using this visualization of $\Lambda$, we observe that an element $X \in \mathsf{mod}\left(\Lambda\otimes_K A_3^{CT}\right)$ can be considered as three $\Lambda$-modules $M_1, M_2, M_3 \in \mathsf{mod} \Lambda$ and a chain of morphisms $M_1 \leftarrow M_2 \leftarrow M_3 \leftarrow M_1$ so that the composition of two consecutive morphisms is zero.
This observation motivates the following definition.

\begin{definition}
	Let
	$$X = \cdots \xleftarrow{d_{0}} X_0 \xleftarrow{d_1} X_1 \xleftarrow{d_2} X_2 \xleftarrow{d_3} \cdots$$
	be a chain complex over $\mathsf{mod} \Lambda$. By definition, $X$ is pointwise finite; that is, each $X_i \in \mathsf{mod} \Lambda$. We say $X$ is \emph{3-periodic} if we have $d_m = d_n$ whenever $m \equiv n \pmod{3}$. Given $X, Y$ two 3-periodic chain complexes, we define a morphism of 3-periodic chain complexes to be a chain map $f: X \rightarrow Y$ such that $f_m = f_n$ whenever $m\equiv n \pmod{3}$. Under these definitions, we denote by $C^3(\Lambda)$ the category of 3-periodic chain complexes over $\mathsf{mod}\Lambda$.
\end{definition}
Based on the discussion above, we conclude:

\begin{proposition}
	The category $\mathsf{mod}\left(\Lambda\otimes_K A_3^{CT}\right)$ is isomorphic to the category $C^3(\Lambda)$ of 3-periodic chain complexes over $\mathsf{mod} \Lambda$.
\end{proposition}

Based on this fact, we will identify $C^3(\Lambda)$ with $\mathsf{mod} \left(\Lambda\otimes_K A_3^{CT}\right)$ for the remainder of this paper. We likewise use $C(\Lambda)$ and $C^b(\Lambda)$ to refer to the categories of (pointwise finite) chain complexes and bounded chain complexes over $\mathsf{mod}\Lambda$. We remark that if we let $A_{\infty}$ be the quiver with vertex set $\mathbb{Z}$ and an arrow $i \rightarrow i-1$ for all $i$, we can also identify $C^b(\Lambda)$ with $\mathsf{mod}\left( \Lambda\otimes_K (KA_\infty/\mathsf{rad}^2 KA_\infty)\right)$ if we choose.

Before proceeding, we set some notation for clarity.

\begin{notation}\
	\begin{enumerate}
		\item Let $i \in \mathbb{Z}$. We denote by $[i]$ the equivalence class of $i\pmod{3}$.
		\item Let $X \in C^3(\Lambda)$. If we write
			$$X = M \leftarrow M' \leftarrow M'' \leftarrow M$$
			without specifying the degree of any of the modules, we assume $M$ is in degree $[-1]$. That is, $M$ is in degree $i$ for all $i \equiv -1\pmod{3}$. We will usually write this as $X_{[-1]} = M$.
		\item Let $X \in C^b(\Lambda)$. If we write
	$$X = M \leftarrow M' \leftarrow \cdots \leftarrow M''$$
	without specifying the degree of any of the modules, we assume $M$ is in degree -1.
	\end{enumerate}
\end{notation}

Our next goal is to relate projective covers and syzygies in $C^b(\Lambda)$ and $C^3(\Lambda)$. We start with a discussion of how to ``wrap'' a bounded chain complex into a 3-periodic chain complex.

\begin{definition}
	Let $X = (M_i, d_i) \in C^b(\Lambda)$. We define a 3-periodic chain complex
		$$WX = \left(\bigoplus_{n\equiv i\pmod{3}} M_n, \bigoplus_{n\equiv i\pmod{3}} d_n\right)$$
		which we call the ``wrapping" of $X$.
\end{definition}

By the assumption that $X$ is bounded, we observe that $WX$ is pointwise finite, so it is a 3-periodic chain complex as desired. Moreover, it is clear that if there is a map $X\xrightarrow{f} Y$ in $C^b(A)$, there is an induced map $WX\xrightarrow{Wf}WY$ given by again taking the direct sum over degrees congruent mod 3. Thus, we have the following.

\begin{lemma}
	Wrapping is a functor $W: C^b(\Lambda) \rightarrow C^3(\Lambda)$. Moreover, this functor is exact.
\end{lemma}

\begin{proof}
	It is clear that $W$ is a functor. To see that $W$ is exact, let
		$$(M_i,d_i) \hookrightarrow (E_i,e_i) \twoheadrightarrow (N_i,f_i)$$
	be an exact sequence of bounded chain complexes. Equivalently, for $i \in \mathbb{Z}$, the induced sequence of modules
		$$M_i \hookrightarrow E_i \twoheadrightarrow N_i$$
		is exact. Thus by the exactness of the bifunctor $\bigoplus$, for each equivalence class $i \in \mathbb{Z}/(3)$, we have an exact sequence
		$$\bigoplus_{n \equiv i\pmod{3}}M_i \hookrightarrow \bigoplus_{n \equiv i\pmod{3}}E_i \twoheadrightarrow \bigoplus_{n \equiv i\pmod{3}}N_i.$$
		Therefore the sequence 
		$$W(M_i,d_i) \hookrightarrow W(E_i,e_i)\twoheadrightarrow W(N_i,f_i)$$
		is exact as well.
\end{proof}

We now give several examples and basic properties of wrapping to better illustrate the concept.

\begin{example}\
	\begin{enumerate}
		\item The wrapping of the bounded chain complex
		$$X = 0 \leftarrow M_0 \xleftarrow{d_0} M_1 \xleftarrow{d_1} M_2 \xleftarrow{d_2} M_3 \xleftarrow{d_3} M_4 \leftarrow 0$$
	in $C^b(\Lambda)$ is the 3-periodic chain complex
		$$WX = M_2 \xleftarrow{\begin{pmatrix}0 & d_2\end{pmatrix}} M_0 \oplus M_3 \xleftarrow{\begin{pmatrix}d_0 & 0\\0 & d_3\end{pmatrix}} M_1 \oplus M_4\xleftarrow{\begin{pmatrix}d_1 \\ 0\end{pmatrix}} M_2$$
		in $C^3(\Lambda)$.
		\item There exist non-isomorphic bounded chain complexes which become isomorphic after wrapping. For example, given any nonzero $M \in \mathsf{mod} \Lambda$, the bounded chain complex consisting only of $M$ in degree 0 and the bounded chain complex consisting only of $M$ in degree 3 become isomorphic after wrapping.
		\item Not every 3-periodic chain complex can be obtained by wrapping a bounded chain complex; that is, $W$ is not essentially surjective. For example, the 3-periodic chain complex
		$$\begin{matrix}1\\2\end{matrix} \leftarrow \begin{matrix}2\\3\end{matrix}\leftarrow\begin{matrix}3\\1\end{matrix}\leftarrow \begin{matrix}1\\2\end{matrix}$$
		in $C^3(A_3^{CT})$, where each boundary map is nontrivial, is not isomorphic to anything in the image of $W$.
		\item Since $\bigoplus$ is faithful, $W$ is as well; however, $W$ is not full. For example, given any nonzero $M \in \mathsf{mod} \Lambda$, the chain complex
		$$X = M \leftarrow 0 \leftarrow 0 \leftarrow M$$
		has $\dim\mathsf{End}_{C^b(\Lambda)}(X) = 2\cdot\dim\mathsf{End}_\Lambda(M)$ and the corresponding
		$$WX = M\oplus M \leftarrow 0 \leftarrow 0 \leftarrow M\oplus M$$
		has $\dim\mathsf{End}_{C^3(\Lambda)}(X) = 4\cdot\dim\mathsf{End}_\Lambda(M)$.
	\end{enumerate}
\end{example}

The following lemma is based on a similar result in \cite[Chapter 3]{happel_triangulated} which appears in Happel's construction of bounded derived categories. The key difference is that Happel considers projectives with respect to the split exact structure, whereas we are considering projectives with respect to the standard exact structure.

\newpage

\begin{lemma}\
	\begin{enumerate}
		\item The indecomposable projective objects of $C^b(\Lambda)$, up to shifting degrees, are precisely objects of the form
	$$0 \leftarrow P \xleftarrow{1} P \leftarrow 0$$
	for $P$ an indecomposable projective in $\mathsf{mod} \Lambda$.
		\item The indecomposable projective objects of $C^3(\Lambda)$, up to cyclic permutation of the degrees, are precisely objects of the form
	$$0 \leftarrow P \xleftarrow{1} P \leftarrow 0$$
	for $P$ an indecomposable projective in $\mathsf{mod} \Lambda$.
	\end{enumerate}
\end{lemma}

We now observe that the functor $W$ sends projectives to projectives. Moreover, we observe that if $0 \neq X \in C^b(\Lambda)$ then $0 \neq WX \in C^3(\Lambda)$. In particular, we have the following.

\begin{lemma}\label{lem:commutes} Let $X \in C^b(\Lambda)$ and let $X \xleftarrow{q} P_X$ be the projective cover of $X$.
	\begin{enumerate}
		\item Then $WX \xleftarrow{Wq} WP_X$ is the projective cover of $WX$ in $C^3(\Lambda)$. That is, taking projective covers commutes with wrapping.
		\item Then $\Omega_{C^3(\Lambda)} WX \cong W\Omega_{C^b(\Lambda)} X$. That is, taking syzygies commutes with wrapping.
	\end{enumerate}
\end{lemma}

\begin{proof}
	(1) Let $X = (M_i,d_i)$ in $C^b(\Lambda)$ with projective cover $(P_X,q)$. As $W$ is exact and preserves projectives, $Wq:WP_X \rightarrow WX$ is a surjection from a projective object to $WX$. The fact that this is a projective cover is clear, because if $WP_X$ had a superfluous direct summand, $P_X$ would as well.
	
	(2) Follows from (1) and the exactness of $W$.
\end{proof}

It is, in general, much more difficult to compute syzygies in $C^3(\Lambda)$ than in $C^b(\Lambda)$. Thus, in light of Lemma \ref{lem:commutes}, we will proceed (in Section \ref{sec:computation}) by finding two families of chain complexes, $\{X_k\}_k$ and $\{Y_k\}_k$, in $C^b(A)$ so that $\Omega^t_{C^b(\Lambda)}(X_k)$ and $\Omega^t_{C^b(\Lambda)}(Y_k)$ first become isomorphic for arbitrarily large values of $t$. We will then pass by wrapping to $C^3(\Lambda)$ and show this is still the case. Before we do this, we more closely examine the computations of certain syzygies in the category $C^b(\Lambda)$.

\section{Syzygies in the Category of Chain Complexes}\label{sec:syzygy}
In this section, we study the syzygies of certain nice chain complexes corresponding to `truncated projective resolutions', which we make precise below. We still assume that $\Lambda$ is an elementary $K$-algebra. For clarity, we denote a morphism $f:\bigoplus_i M_i \rightarrow \bigoplus_j N_j$ by a set of arrows $\{f_{i,j}: M_i \rightarrow N_j\mid N_j \cap \mathsf{image}(f|_{M_i}) \neq 0\}$. We will also omit the quotient maps from the data of projective covers.

\begin{definition}
	Choose some integer $m > 0$ and let $M \in \mathsf{mod}\Lambda$ be of projective dimension at least $m$. Let $P_0$ be the projective cover of $M$. Now let $Q_1 \oplus R_1$ be a direct sum decomposition of $\ker(P_0 \rightarrow M)$ so that $R_1 \neq 0$ has no projective direct summand and $\mathsf{pd}_{\Lambda}(R_1) \geq m-1$. Let $P_1$ be the projective cover of $R_1$. Likewise, for $1 < i < m$, define $P_i$ and $Q_i$ inductively starting with $\ker(P_{i-1}\rightarrow R_{i-1})$. Denote $Q_m = \ker(P_{m-1}\rightarrow R_{m-1})$. This data, arranged in the chain complex shown below, gives an object $X \in C^b(\Lambda)$ which we call a \emph{truncated (minimal) projective resolution} of $M$ (of length $m$).
\end{definition}
	\begin{center}
	\begin{tikzcd}[column sep = 7mm,row sep = 2mm]
		&&Q_1 \arrow[d,no head,white,anchor=center,"\oplus"{black,description}]& Q_2 \arrow[d,no head,white,anchor=center,"\oplus"{black,description}]& \cdots & Q_{m-1} \arrow[d,no head,white,anchor=center,"\oplus"{black,description}]& Q_m \\
		M & P_0 \arrow[l,two heads] \arrow[ur,hookleftarrow] & P_1 \arrow[l] \arrow[ur,hookleftarrow] & P_2 \arrow[l]&\arrow[ur,hookleftarrow] \cdots \arrow[l] & \arrow[ur,hookleftarrow] P_{m-1}\arrow[l]\\
	\end{tikzcd}
	\end{center}
	
Throughout this paper, every truncated projective resolution will be minimal. In this section we compute the syzygy of an arbitrary truncated projective resolution. In Section \ref{sec:computation} we pass to the special case of our counterexample. We first prove the following, which will be critical in Section \ref{sec:computation}.

\begin{lemma}\label{lem:indecomposable}
	Let $X = (X_i,d_i)$ be a bounded chain complex. By shifting as necessary, suppose without loss of generality that $X_i = 0$ for $i < -1$. Suppose further that
	\begin{enumerate}
		\item $X_{-1}$ is an indecomposable $\Lambda$-module.
		\item For all $i\equiv j\pmod{3}$ with $i \neq -1$, we have $\mathsf{Hom}_\Lambda(X_i,\ker(d_j)) = 0$.
	\end{enumerate}
	Then $WX$ is indecomposable in $C^3(\Lambda)$.
\end{lemma}

\begin{proof}
	Let $m$ be the largest index for which $X_m \neq 0$. We prove the result for $m\equiv0\pmod{3}$, although the other cases follow analogously. We first observe that we can write
	\begin{center}
	\begin{tikzcd}[row sep = 1.5mm]
		&X_{-1} & X_0\arrow[l, two heads,"d_0" above] & X_1\arrow[l, "d_1" above] & X_{-1}\\
		\arrow[r,no head,white,anchor=center,"\textnormal{\normalsize $WX=$}"{black,description},yshift=-3mm,xshift = -12mm]& X_2& X_3\arrow[l,"d_3" above] & X_4 \arrow[l,"d_4" above] & X_2\arrow[ul,"d_5" above]\\
		&\vdots&\vdots&\vdots&\vdots\arrow[ul,"d_8"above]\\
		&X_{m-1}& X_m\arrow[l,left hook->,"d_{m}"above] && X_{m-1}\arrow[ul,"d_{m-1}"above,pos=.4]
	\end{tikzcd}
	\end{center}
	For $-1 \leq \ell \leq m$, we denote by $\pi_\ell$ the projection $WX_{[\ell]} \twoheadrightarrow X_\ell$. Now suppose for a contradiction that there is a direct sum decomposition $WX \cong B \oplus C$. Then without loss of generality, (1) implies that $C_{[-1]}$ does not contain any direct summand of $X_{-1}$ as a direct summand. We claim that this implies $C = 0$ and hence $WX \cong B$ is indecomposable.
	
	Let $\iota = (\iota_{[-1]},\iota_{[0]},\iota_{[1]}): C \hookrightarrow X$ be the inclusion map. We will show that $\iota = 0$ by showing $\pi_{\ell}\circ \iota_{[\ell]} = 0$ for all $-1 \leq \ell \leq m$. First observe that by taking $j = -1$, (2) implies $\mathsf{Hom}(C_{[-1]},X_{-1}) = 0$ and thus $\pi_{-1}\circ \iota_{[-1]} = 0$.
	
	Assume by induction on $\ell$ that the result holds for $-1 \leq \ell < m$. We then have a commutative diagram
	\begin{center}
		\begin{tikzcd}[row sep = small]
			C_{[\ell]}\arrow[d,"\iota_{[\ell]}" left] & C_{[\ell+1]} \arrow[l,"e_{[\ell+1]}"above]\arrow[d,"\iota_{[\ell + 1]}" right]\\
			X_{[\ell]} \arrow[d,"\pi_\ell" left]& X_{[\ell+1]}\arrow[d,"\pi_{\ell+1}"right]\arrow[l,"d_{[\ell+1]}" above]\\
			X_\ell & \arrow[l,"d_{\ell+1}" below] X_{\ell + 1}
		\end{tikzcd}
	\end{center}
	where $e_{[\ell+1]}$ is the boundary map of $C$ and $\displaystyle d_{[\ell + 1]} = \bigoplus_{j \equiv \ell+1\pmod{3}} d_j$. By assumption, we have $\pi_{\ell} \circ \iota_{[\ell]} \circ e_{[\ell+1]} = 0$. Since $\mathsf{Hom}_A(C_{[\ell+1]},\ker(d_{\ell+1})) = 0$ by (2), this implies that $\pi_{\ell + 1}\circ \iota_{[\ell + 1]} = 0$ as desired. We conclude that $\iota = 0$ and hence $C = 0$.
\end{proof}

We now prove the main result of this section.

\begin{proposition}\label{prop:syzygy}
	Let $M \in \mathsf{mod} \Lambda$ be of projective dimension at least $m$ and let $X$ be a truncated projective resolution of $M$ of length $m$, as shown below. For each $i \in \{1,\ldots,m\}$, let $P_{Q_i}$ be the projective cover of $Q_i$. Then $\Omega_{C^b(\Lambda)} X$ and the projective cover of $X$ are as shown below.
	
	\begin{center}
	\begin{tikzcd}[column sep = 4mm,row sep = 1.5mm]
	
		\phantom{P}\arrow[dd,start anchor = north, end anchor = south, no head, xshift=-5mm,decorate, decoration={brace,mirror},"\textnormal{\large $\Omega_{C^b(\Lambda)} X = $}" left=1mm]&& \Omega_\Lambda Q_1\arrow[d,no head,white,anchor=center,"\oplus"{black,description}]& \Omega_\Lambda Q_2\arrow[d,no head,white,anchor=center,"\oplus"{black,description}] & \cdots & \Omega_\Lambda Q_{m-1}\arrow[d,no head,white,anchor=center,"\oplus"{black,description}] & \Omega_\Lambda Q_m\arrow[dd,start anchor = north, end anchor = south, no head, xshift=5mm,decorate, decoration={brace}]\\
		\phantom{P} & P_{Q_1} \arrow[ur,hookleftarrow,"j_1"sloped,pos=1] \arrow[d,no head,white,anchor=center,"\oplus"{black,description}]& P_{Q_2}\arrow[d,no head,white,anchor=center,"\oplus"{black,description}] \arrow[ur,hookleftarrow,"j_2"sloped,pos=0.9]& P_{Q_3}\arrow[d,no head,white,anchor=center,"\oplus"{black,description}]\arrow[ur,hookleftarrow,"j_3"sloped,pos=0.8]&\cdots\arrow[ur,hookleftarrow,"j_{m\text-1}"sloped,pos=1]&P_{Q_{m}}\arrow[ur,hookleftarrow,"j_m"sloped,pos=0.8] & \phantom{P}\\
		\Omega_A M \arrow[ur,leftarrow,"i_1\circ q_1" sloped,pos=1.4] & P_1\arrow[l,"f_1"]\arrow[ur,leftarrow,"i_2\circ q_2"sloped,pos=1.1] & P_2 \arrow[l,"f_2"] \arrow[ur,leftarrow,"i_3\circ q_3"sloped,pos=1]& P_3\arrow[l,"f_3"]\arrow[ur,leftarrow,"i_4\circ q_4"sloped,pos=1.1]& \cdots\arrow[l,"f_4"]\arrow[ur,leftarrow,"i_m\circ q_m"sloped,pos=1.2] & \phantom{P} & \phantom{P}\\[0.75em]
		&&&&\cdots&&\\[0.75em]
	
		\phantom{P}\arrow[ddd,start anchor = north, end anchor = south, no head, xshift=-5mm,decorate, decoration={brace,mirror},"\textnormal{\large $P_X=$}" left=5mm]\arrow[uu,hookleftarrow,"h_{-1}" right] &\phantom{P}\arrow[uu,hookleftarrow,"h_{0}" right] & P_{Q_1}\arrow[d,no head,white,anchor=center,"\oplus"{black,description}]\arrow[uu,hookleftarrow,"h_{1}" right]  & P_{Q_2}\arrow[d,no head,white,anchor=center,"\oplus"{black,description}] \arrow[uu,hookleftarrow,"h_{2}" right] & \cdots & P_{Q_{m-1}}\arrow[d,no head,white,anchor=center,"\oplus"{black,description}] \arrow[uu,hookleftarrow,"h_{m-1}" right] & P_{Q_m}\arrow[uu,hookleftarrow,"h_{m}" right] \arrow[ddd,start anchor = north, end anchor = south, no head, xshift=5mm,decorate, decoration={brace}]\\
		\phantom{P}\arrow[uuuu,hookleftarrow,xshift=-17mm,"\textnormal{\normalsize $h$}"left=1mm]& P_0\arrow[d,no head,white,anchor=center,"\oplus"{black,description}] & P_1\arrow[d,no head,white,anchor=center,"\oplus"{black,description}] & P_2\arrow[d,no head,white,anchor=center,"\oplus"{black,description}] & \cdots & P_{m-1}\arrow[d,no head,white,anchor=center,"\oplus"{black,description}]\\
		\phantom{P}\arrow[ddd,two heads,xshift=-17mm,"\textnormal{\normalsize$g$}"left=1mm]&P_{Q_1}\arrow[d,no head,white,anchor=center,"\oplus"{black,description}]  \arrow[uur,leftarrow,"1",pos=.6] & P_{Q_2}  \arrow[d,no head,white,anchor=center,"\oplus"{black,description}]\arrow[uur,leftarrow,"1",pos=.6] & P_{Q_3}  \arrow[d,no head,white,anchor=center,"\oplus"{black,description}]\arrow[uur,leftarrow,"1",pos=.6] & \cdots  \arrow[uur,leftarrow,"1",pos=.6] &  P_{Q_m} \arrow[uur,leftarrow,"1",pos=.6]\\
		P_0  \arrow[uur,leftarrow,"1"] \arrow[dd,two heads,"g_{-1}"]& P_1  \arrow[uur,leftarrow,"1",pos=.6] \arrow[dd,two heads,"g_0"]& P_2  \arrow[uur,leftarrow,"1",pos=.6] \arrow[dd,two heads,"g_1"]& P_3  \arrow[uur,leftarrow,"1",pos=.6] \arrow[dd,two heads, "g_2"]& \cdots   & \phantom{P} \arrow[dd,two heads,"g_{m-1}"] & \phantom{P} \arrow[dd,two heads,"g_m"]\\[0.75em]
		&&&&\cdots\\[0.75em]
		
		\phantom{P}\arrow[d,start anchor = north, end anchor = south, no head, xshift=-5mm,decorate, decoration={brace,mirror},"\textnormal{\large $X=$}" left=5.5mm]&\phantom{P}&Q_1\arrow[d,no head,white,anchor=center,"\oplus"{black,description}] & Q_2 \arrow[d,no head,white,anchor=center,"\oplus"{black,description}]& \cdots & Q_{m-1} \arrow[d,no head,white,anchor=center,"\oplus"{black,description}]& Q_m\arrow[d,start anchor = north, end anchor = south, no head, xshift=5mm,decorate, decoration={brace}]\\
		M & P_0 \arrow[l,two heads,"f_0"] \arrow[ur,hookleftarrow,"i_1"sloped,pos=.7] & P_1 \arrow[l,"f_1"]\arrow[ur,hookleftarrow,"i_2"sloped,pos=.7] & P_2 \arrow[l,"f_2"]\arrow[ur,hookleftarrow,"i_3"sloped,pos=.7] &\arrow[ur,hookleftarrow,"i_{m-1}"sloped,pos=1.1] \cdots \arrow[l,"f_3"]& \arrow[ur,hookleftarrow,"i_m"sloped,pos=.8] P_{m-1}\arrow[l,"f_{m-1}"]&\phantom{P}\\
	\end{tikzcd}
	\end{center}
	For all $k$, the maps $i_k$ and $j_k$ are inclusion maps and the map $q_k:P_{Q_k} \rightarrow Q_k$ is the quotient map. Likewise, $\iota: \Omega_A(M)\hookrightarrow P_0$ is the inclusion map. Moreover,
		$$h_k = \begin{cases}
		-\iota & k = -1\\
		\begin{pmatrix}-i_1\circ q_1 & -f_1\\1 & 0\\0 & 1\end{pmatrix} & k = 0\\
		(-1)^k\cdot\begin{pmatrix} -j_k & 0 & 0\\0 & -i_{k+1}\circ q_{k+1} & -f_{k+1}\\0 & 1 & 0\\0 & 0 & 1\end{pmatrix} & 1 \leq k \leq m-2\\
		(-1)^{m-1}\cdot\begin{pmatrix} -j_{m-1} & 0\\0 & -i_{m}\circ q_{m}\\0 & 1\end{pmatrix} & k = m-1\\
		(-1)^{m+1} j_m & k = m
		\end{cases}	
	$$and
	$$g_k = \begin{cases}
		f_0 & k = -1\\
		\begin{pmatrix}1 & i_1\circ q_1 & f_1\end{pmatrix} & k = 0\\
		\begin{pmatrix} q_k & 0 & 0 & 0\\0 & 1 & i_{k+1}\circ q_{k+1} & f_k\end{pmatrix} & 1 \leq k \leq m-2\\
		\begin{pmatrix}q_{m-1} & 0 & 0\\0 & 1 & i_{m}\circ q_{m}\end{pmatrix} & k = m-1\\
		q_m & k = m
		\end{cases}$$
\end{proposition}

\begin{proof}
	It can be verified directly that the diagram commutes and every column corresponds to an exact sequence.
\end{proof}

As a direct consequence of Proposition \ref{prop:syzygy}, we have the following.

\newpage

\begin{corollary}\label{cor:formulaResults} Let $M$ and $X$ be as above.
	\begin{enumerate}
		\item $\Omega_{C^b(\Lambda)}X$ is a truncated projective resolution of $\Omega_\Lambda M$. It is of length $m-1$ if $Q_m$ is projective and is of length $m$ otherwise. In particular, higher syzygies of $X$ can be computed using the procedure outlined in the proposition.
		\item Suppose $Q_1 \neq 0$. Then $Q_1 \leftarrow P_{Q_1} \leftarrow \Omega_A Q_1$ is a direct summand of $\Omega_{C^b(\Lambda)} X$.
		\item Suppose $Q_k \neq 0$. Then there is a truncated projective resolution of the form
		$$Q_k \leftarrow P'_0 \leftarrow P'_1 \leftarrow \cdots \leftarrow P'_{k-1}\leftarrow \Omega_A^k Q_k$$
		which is a direct summand of $\Omega_{C^b(\Lambda)}^k X$.
	\end{enumerate}
\end{corollary}

In the next section, we will use these results to study the syzygies of specific truncated projective resolutions related to our counterexample.


\section{Computations for the Counterexample}\label{sec:computation}
We now fix $A = KQ/\mathsf{rad}^2KQ$ where $Q$ is the quiver
\begin{center}
\begin{tikzcd}[row sep = tiny]
		&2\arrow[dl,"x_2" above left]\arrow[dr,"x_2'"above right]\\
		3\arrow[dr,"x_3"below left]&&4\arrow[dl,"x_4"below right]&\\
		&1\arrow[uu,"x_1"right]
\end{tikzcd}
\end{center}
as in the introduction. The goal of this section is to prove our main theorem:

\begin{theorem}\label{thm:phidiminfinite}
	The $\phi$-dimension of $A\otimes_K A_3^{CT}$ is infinite.
\end{theorem}

Throughout this section, we will omit direct sum symbols in our depictions of chain complexes. As is standard, for $i \in \{1,2,3,4\}$, we denote by $S_i,P_i \in \mathsf{mod} A$ the simple module supported at vertex $i$ and its projective cover. Each of these eight modules is indecomposable and has endomorphism ring isomorphic to the ground field $K$. With this notation, the minimal projective resolution of the simple module $S_3 \in \mathsf{mod} A$ is as follows, and the minimal projective resolution of $S_4$ is the same with all of the indices 3 and 4 interchanged. 
\begin{center}
\begin{tikzcd}[row sep = 1.5mm,column sep = 5mm]
	&&&&&&& P_3 \arrow[dl] & \cdots \arrow[l]\\
	&&&& P_3 \arrow[dl] & P_1 \arrow[l] & P_2 \arrow[l] & P_4\arrow[l] & \cdots\arrow[l]\\
	S_3 & P_3 \arrow[l,two heads] & P_1\arrow[l] & P_2 \arrow[l] & P_4 \arrow[l] & P_1 \arrow[l] & P_2 \arrow[l] & P_3 \arrow[l] & \cdots \arrow[l]\\
	&&&&&&& P_4 \arrow[ul] & \cdots \arrow[l]\\
\end{tikzcd}
\end{center}
 We begin with an example to motivate this process. We remark that 
 	$$\Omega_A S_i \cong \begin{cases} S_2 & i = 1\\S_3\oplus S_4 & i = 2\\S_1 & i = 3, 4.\end{cases}$$
In particular, all four simple modules are of infinite projective dimension.

\begin{example}\label{ex:ell=1}
	Let $X_1$ be the truncated projective resolution of $S_3$ given below.
	\begin{center}
	\begin{tikzcd}[row sep = 1.5mm,column sep = 5mm]
		&&&&& S_3 \arrow[dl, left hook->]&&&S_3\arrow[dl, left hook->]\\
		X_1 = &S_3 & P_3 \arrow[l,two heads] & P_1\arrow[l] & P_2 \arrow[l] & P_4 \arrow[l] & P_1 \arrow[l] & P_2 \arrow[l] & S_4 \arrow[l,left hook->]\\
	\end{tikzcd}
	\end{center}
	Let $Y_1$ be obtained from $X_1$ by interchanging all indices 3 and 4.
	By Proposition \ref{prop:syzygy}, we have
	\begin{center}
	\begin{tikzcd}[row sep = 1.5mm,column sep = 5mm]
		&&&& P_3 \arrow[dl]&S_1\arrow[l,left hook->]&&P_3 \arrow[dl] &S_1\arrow[l,left hook->]\\
		\Omega_{C^b(A)} X_1 = &S_1 & P_1 \arrow[l,two heads] & P_2\arrow[l] & P_4 \arrow[l] & P_1 \arrow[l] & P_2 \arrow[l] & P_4 \arrow[l] & S_1 \arrow[l, left hook->]\\
	\end{tikzcd}
	\end{center}
	Applying the proposition two more times, we see that
	\begin{center}
	\begin{tikzcd}[row sep = 1.5mm,column sep = 5mm]
		&&&&& S_3 \arrow[dl, left hook ->]&&&S_3 \arrow[dl,left hook->]\\
		\arrow[r,no head,white,anchor=center,"\textnormal{\normalsize $\Omega_{C^b(A)}^3 X_1 =$}"{black,description},yshift=-4mm,xshift = -12mm]&S_3 & P_3 \arrow[l,two heads] & P_1 \arrow[l] & P_2 \arrow[l] & S_4 \arrow[l,left hook->] &P_1\arrow[dl]&P_2\arrow[l] &S_4 \arrow[l,left hook->]\\
		&S_4 & P_4 \arrow[l,two heads] & P_1\arrow[l] & P_2 \arrow[l] & P_3 \arrow[l] && & S_3 \arrow[dl, left hook->]\\
		&&&&&P_4\arrow[ul]&P_1\arrow[l]&P_2\arrow[l]&S_4 \arrow[l, left hook->]\\
	\end{tikzcd}
	\end{center}
	In particular, this breaks apart into a direct sum of two chain complexes. Now certainly, we see that $\Omega_{C^b(A)}^3 X_1 \ncong\Omega_{C^b(A)}^3 Y_1$. We will show in Proposition \ref{prop:phidiminfinite} this is still the case after we apply the functor $W$; that is, $\Omega_{C^3(A)}^3 WX_1 \ncong \Omega_{C^3(A)}^3 WY_1$. Moreover, applying Proposition \ref{prop:syzygy} once more, we obtain
	\begin{center}
	\begin{tikzcd}[row sep = 1.5mm,column sep = 5mm]
		&&&&P_3\arrow[dl]& S_1 \arrow[l, left hook ->]&&P_3\arrow[dl]&S_1 \arrow[l,left hook->]\\
		\arrow[r,no head,white,anchor=center,"\textnormal{\normalsize $\Omega_{C^b(A)}^4 X_1 =$}"{black,description},yshift=-4mm,xshift = -12mm]&S_1 & P_1 \arrow[l,two heads] & P_2 \arrow[l] & P_4 \arrow[l] & S_1 \arrow[l,left hook->] &P_2\arrow[dl]&P_4\arrow[l] &S_1 \arrow[l,left hook->]\\
		&S_1 & P_1 \arrow[l,two heads] & P_2\arrow[l] & P_3 \arrow[l] & P_1 \arrow[l] && P_3\arrow[dl]& S_1 \arrow[l, left hook->]\\
		&&&&P_4\arrow[ul]&P_1\arrow[l]&P_2\arrow[l]&P_4\arrow[l]&S_1 \arrow[l, left hook->]\\
	\end{tikzcd}
	\end{center}
	We observe that $\Omega_{C^b(A)}^4 X_1 \cong \Omega_{C^b(A)}^4 Y_1$ because both complexes are symmetric in the indices 3 and 4. In particular, these complexes will still be isomorphic when we pass via wrapping to the corresponding 3-periodic chain complexes.
\end{example}

We observe that the rightmost asymmetry (between the indices 3 and 4) of $X_1$ occurs in degree 3 and that this asymmetry vanishes after precisely 4 applications of $\Omega_{C^b(A)}$. Definition \ref{def:XandY} and Proposition \ref{prop:phidiminfinite} below generalize this observation to longer truncated projective resolutions. The key to Proposition \ref{prop:phidiminfinite} is that this pattern persists after application of the wrapping functor, giving a systematic way to construct 3-periodic chain complexes of arbitrarily high $\phi$-dimension.

\begin{definition}\label{def:XandY}Let $k \geq 0$.
\begin{enumerate}
	\item We define:
	\begin{center}
	\begin{tikzcd}[row sep = 1mm, column sep = 5mm]
		&&&\phantom{P}\arrow[d,start anchor = north, end anchor = south, no head, xshift=-2mm,decorate, decoration={brace,mirror}]&& S_3\arrow[d,start anchor = north, end anchor = south, no head, xshift=2mm,decorate, decoration={brace},"\textnormal{\normalsize$k$}"{yshift=4mm,xshift = 1mm}] \arrow[dl, left hook->]&&&S_3\arrow[dl, left hook->]\\
		X_k = &S_3 & P_3 \arrow[l,two heads] & P_1\arrow[l] & P_2 \arrow[l] & P_4 \arrow[l] & P_1 \arrow[l] & P_2 \arrow[l] & S_4 \arrow[l,left hook->]\\
	\end{tikzcd}
	\begin{tikzcd}[row sep = 1mm, column sep = 5mm]
		&&&\phantom{P}\arrow[d,start anchor = north, end anchor = south, no head, xshift=-2mm,decorate, decoration={brace,mirror}]&& S_4\arrow[d,start anchor = north, end anchor = south, no head, xshift=2mm,decorate, decoration={brace},"\textnormal{\normalsize$k$}"{yshift=4mm,xshift = 1mm}] \arrow[dl, left hook->]&&&S_4\arrow[dl, left hook->]\\
		Y_k = &S_4 & P_4 \arrow[l,two heads] & P_1\arrow[l] & P_2 \arrow[l] & P_3 \arrow[l] & P_1 \arrow[l] & P_2 \arrow[l] & S_3 \arrow[l,left hook->]\\
	\end{tikzcd}
	\end{center}
	The notation $\{-\}^k$ means the portion in the brackets is repeated $k$ times connected by arrows. For example,
	\begin{center}
	\begin{tikzcd}[row sep = 1.5mm,column sep = 2.5mm]
		&&&\phantom{P}\arrow[d,start anchor = north, end anchor = south, no head, xshift=-2mm,decorate, decoration={brace,mirror}]&& S_3\arrow[dl, left hook->]\arrow[d,start anchor = north, end anchor = south, no head, xshift=2mm,decorate, decoration={brace}] &\phantom{P}\arrow[d,start anchor = north, end anchor = south, no head, xshift=-2mm,decorate, decoration={brace,mirror}]&&S_3\arrow[dl, left hook->]\arrow[d,start anchor = north, end anchor = south, no head, xshift=2mm,decorate, decoration={brace}] &&&S_3\arrow[dl,left hook->]\\
		X_2 = &S_3 & P_3 \arrow[l,two heads] & P_1\arrow[l] & P_2 \arrow[l] & P_4 \arrow[l] & P_1 \arrow[l] & P_2 \arrow[l] & P_4 \arrow[l] & P_1 \arrow[l] & P_2 \arrow[l] &S_4 \arrow[l, left hook->]\\
	\end{tikzcd}
	\end{center}
	
	\item Let $i \in \{1,2,3,4\}$. We define $Z_k^i$ to be the chain complex formed by truncating the projective resolution of $S_i$ at degree $3+3k$. For example,
	\begin{center}
		\begin{tikzcd}[row sep = 1.5mm, column sep = 5mm]
				&&&&&&&&S_3\arrow[dl,left hook->]\\
				&&&&&P_3 \arrow[dl] &P_1\arrow[l]&P_2\arrow[l] &S_4 \arrow[l,left hook->]\\
				Z_1^4 = &S_4 & P_4 \arrow[l,two heads] & P_1\arrow[l] & P_2 \arrow[l] &&& & S_3 \arrow[dl, left hook->]\\
				&&&&&P_4\arrow[ul]&P_1\arrow[l]&P_2\arrow[l]&S_4 \arrow[l, left hook->]\\
		\end{tikzcd}
	\end{center}
\end{enumerate}
\end{definition}

We remark that we could have chosen to define these `$Z$-type' complexes to be truncated at any degree rather than just at degrees of the form $3 + 3k$. Truncation at degrees $3+3k$ is for notational convenience in the statement of Lemma \ref{lem:bigSyzygy} below.

We claim the following about these chain complexes:

\begin{proposition}\label{prop:phidiminfinite} Let $k \geq 1$. Then
\begin{enumerate}
	\item $\Omega_{C^3(A)}^{3k+1}WX_k \cong \Omega_{C^3(A)}^{3k+1}WY_k.$
	\item $\Omega_{C^3(A)}^{3k}WX_k \ncong \Omega_{C^3(A)}^{3k} WY_k.$
\end{enumerate}
\end{proposition}
Before proving this, we observe that our main result is an immediate corollary.

\begin{proof}[Proof that Proposition \ref{prop:phidiminfinite} implies Theorem \ref{thm:phidiminfinite}]
	For all $k \in \mathbb{N}$, we have that$$\phi(WX_k\oplus WY_k) \geq 3k$$ by the proposition and Lemma \ref{lem:bound}(2). Therefore
$$\phi\mathsf{dim}(A\otimes_K A_3^{CT}) = \sup\left\{\phi(X)\mid X \in \mathsf{mod} \left(A\otimes_K A_3^{CT}\right)\right\} = \infty$$
as claimed.
\end{proof}

The remainder of this section is devoted to proving Proposition \ref{prop:phidiminfinite}. We begin with several lemmas about these particular chain complexes and their syzygies.

\begin{lemma}\label{lem:indecomposable2}\
	\begin{enumerate}
		\item Let $k \geq 1$. Then $WX_k$ and $WY_k$ are indecomposable in $C^3(A)$.
		\item Let $k \geq 0$ and $i \in \{1,2,3,4\}$. Then $WZ_k^i$ is indecomposable in $C^3(A)$.
	\end{enumerate}
\end{lemma}

\begin{proof}
	Each of these truncated projective resolutions satisfies the hypotheses of Lemma \ref{lem:indecomposable}. For example, let $(M_i,d_i) = Z_0^3$, so we have:
	\begin{center}
		\begin{tabular}{|c||c|c|c|c|c|c}
			\hline
			$i$ & $-1$ & 0 & 1 & 2 & 3\\
			\hline
			$M_i$ & $S_3$ & $P_3$ & $P_1$ & $P_2$ & $S_3\oplus S_4$\\
			\hline
			$\ker(d_i)$ & $S_3$ & $S_1$ & $S_2$ & $S_3 \oplus S_4$ & 0\\
			\hline
		\end{tabular}
	\end{center}
	and all other $M_i$ are trivial. Thus we see that $M_{-1} \cong S_3$ is indecomposable as an $A$-module and for $i \equiv j\pmod{3}$ with $i \neq -1$, there are no nonzero morphisms $M_i \rightarrow \ker(d_j)$, as desired.
\end{proof}

\begin{lemma}\label{lem:smallSyzygy}
	Let $k \geq 0$. Then
		$$\Omega_{C^b(A)} Z_k^i \cong \begin{cases}Z_k^2 & i = 1\\Z_k^3\oplus Z_k^4 & i = 2\\Z_k^1 & i = 3,4.\end{cases}$$
\end{lemma}
\begin{proof}
	This follows immediately from applying the construction in Proposition \ref{prop:syzygy}. For example, if $k = 1$ and $i = 2$, we have
	\begin{center}
		\begin{tikzcd}[row sep = 1mm, column sep = 5mm]
				&&&&&&P_3\arrow[dl] & P_1 \arrow[l] & S_2 \arrow[l,left hook->]\\
				&&&P_3 \arrow[dl] &P_1\arrow[l]&P_2\arrow[l] &P_4 \arrow[l] & P_1 \arrow[l] & S_2 \arrow[l,left hook->]\\
				Z_1^2 = &S_2 & P_2 \arrow[l] &&&& P_3 \arrow[dl] & P_1\arrow[l] & S_2 \arrow[l,left hook->]\\
				&&&P_4\arrow[ul]&P_1\arrow[l]&P_2\arrow[l]&P_4 \arrow[l]&P_1\arrow[l] & S_2\arrow[l,left hook->]\\
		\end{tikzcd}
		\begin{tikzcd}[row sep = 1mm, column sep = 5mm]
				&&&&&P_3 \arrow[dl] &P_1\arrow[l]&P_2\arrow[l] &S_3\oplus S_4 \arrow[l,left hook->]\\
				\arrow[r,no head,white,anchor=center,"\textnormal{\normalsize $\Omega_{C^b(A)}Z_1^2 =$}"{black,description},yshift=-3mm,xshift = -12mm]
				&S_3\arrow[d,start anchor = north, end anchor = south, no head, xshift=-2mm,decorate, decoration={brace,mirror}]& P_3 \arrow[l,two heads] & P_1\arrow[l] & P_2 \arrow[l] &P_4 \arrow[l] & P_1\arrow[l] & P_2\arrow[l] & S_3\oplus S_4 \arrow[l,left hook->]\\
				&S_4 & P_4 \arrow[l,two heads] & P_1\arrow[l] & P_2 \arrow[l] &P_3 \arrow[l] & P_1\arrow[l] & P_2\arrow[l] & S_3\oplus S_4 \arrow[l,left hook->]\\
				&&&&&P_4 \arrow[ul] &P_1\arrow[l]&P_2\arrow[l] &S_3\oplus S_4 \arrow[l,left hook->]\\
		\end{tikzcd}
	\end{center}
	since $\Omega_{A}S_2 \cong S_3 \oplus S_4$.
\end{proof}

\begin{lemma}\label{lem:bigSyzygy}
	Let $k \geq 1$. Then
	\begin{eqnarray*}
		\Omega_{C^b(A)}^{3k}X_k &\cong& Z_k^4 \oplus Z_{k-1}^3 \oplus \left(\displaystyle\bigoplus_{j = 0}^{k-2} (Z_j^3\oplus Z_j^4)^{2^{k-j-2}}\right)\\
		\Omega_{C^b(A)}^{3k}Y_k &\cong& Z_k^3 \oplus Z_{k-1}^4 \oplus \left(\displaystyle\bigoplus_{j = 0}^{k-2} (Z_j^3\oplus Z_j^4)^{2^{k-j-2}}\right).
	\end{eqnarray*}
\end{lemma}

\begin{proof}
	We only show the result for $X_k$, as the proof for $Y_k$ is similar. Recall that 
	\begin{center}
	\begin{tikzcd}[row sep = 1.5mm, column sep = 5mm]
		&&&\phantom{P}\arrow[d,start anchor = north, end anchor = south, no head, xshift=-2mm,decorate, decoration={brace,mirror}]&& S_3\arrow[d,start anchor = north, end anchor = south, no head, xshift=2mm,decorate, decoration={brace},"\textnormal{\normalsize$k$}"{yshift=4mm,xshift = 1mm}] \arrow[dl, left hook->]&&&S_3\arrow[dl, left hook->]\\
		X_k = &S_3 & P_3 \arrow[l,two heads] & P_1\arrow[l] & P_2 \arrow[l] & P_4 \arrow[l] & P_1 \arrow[l] & P_2 \arrow[l] & S_4 \arrow[l,left hook->]\\
	\end{tikzcd}
	\end{center}
	By applying the construction in Proposition \ref{prop:syzygy}, we then have
	\begin{center}
	\begin{tikzcd}[row sep = 1mm,column sep = 2mm]
		&&&&& S_3 \arrow[dl, left hook ->]\\
		&S_3 & P_3 \arrow[l,two heads] & P_1 \arrow[l] & P_2 \arrow[l] & S_4 \arrow[l,left hook->]&&&S_3\arrow[ddd,start anchor = north, end anchor = south, no head, xshift=2mm,decorate, decoration={brace},"\textnormal{\normalsize$k-1$}"{yshift=12mm,xshift = 1mm}] \arrow[dl,left hook->]\\
		\Omega_{C^b(A)}^3 X_k = &&&&&P_3 \arrow[dl]&P_1\arrow[dd,start anchor = north, end anchor = south, no head, xshift=-2mm,decorate, decoration={brace,mirror}]\arrow[l]&P_2\arrow[l] &S_4 \arrow[l,left hook->]\\
		&S_4 & P_4 \arrow[l,two heads] & P_1\arrow[l] & P_2 \arrow[l] & && & P_3 \arrow[dl]&P_1\arrow[l] & P_2\arrow[l] & S_3\oplus S_4 \arrow[l,left hook->]\\
		&&&&&P_4\arrow[ul]&P_1\arrow[l]&P_2\arrow[l]&P_4 \arrow[l]&P_1\arrow[l] & P_2\arrow[l] & S_3\oplus S_4 \arrow[l,left hook->]\\
	\end{tikzcd}
	\end{center}

We first observe that $Z_0^3$ appears as a direct summand of $\Omega^3_{C^b(A)} X_k$. Moreover, by Lemma \ref{lem:smallSyzygy}, we observe that $\Omega^{3(k-1)}Z_0^3 \cong (Z_0^3 \oplus Z_0^4)^{2^{k-2}}$. This accounts for the $j=0$ term in the proposed direct sum decomposition of $\Omega^{3k}_{C^b(A)} X_k$.

Now let $X^1$ be the complement of $Z_0^3$ in $\Omega_{C^b(A)}^3 X_k$. By similar reasoning, we see that $\Omega^3_{C_b(A)} X^1$ will contain $Z_1^3$ as a direct summand. Again applying Lemma \ref{lem:smallSyzygy}, we then have that $\Omega^{3(k-2)}Z_1^3 \cong (Z_1^3 \oplus Z_1^4)^{2^{k-3}}$ is a direct summand of $\Omega^{3k}_{C^b(A)} X_k$.

For $\ell \leq k$, inductively define $X^\ell$ to be the complement of $Z_{\ell-1}^{3}$ in $\Omega^3_{C^b(A)} X^{\ell-1}$. For readability, denote $S_{34}:=S_3\oplus S_4$. Then we see that the general form of $X^{\ell}$ is as shown below.

\begin{center}
\begin{tikzcd}[row sep = 0.5mm, column sep = 2.75mm]
	&&&&&&&&P_3\arrow[dl] & P_1\arrow[ddddddddddddd,start anchor = north, end anchor = south, no head, xshift=-2mm,decorate, decoration={brace,mirror}]\arrow[l] & P_2 \arrow[l] & S_{34}\arrow[l,left hook->]\arrow[ddddddddddddd,start anchor = north, end anchor = south, no head, xshift=2.5mm,decorate, decoration={brace},"\textnormal{\normalsize$k-\ell$}"{yshift=42mm,xshift = 1mm}] \\
	&&&&&&\cdots&P_2\arrow[l]&P_4\arrow[l]& P_1\arrow[l] & P_2 \arrow[l] & S_{34}\arrow[l,left hook->]\\
	&&&&&&\vdots&\vdots&\vdots&\vdots&\vdots&\vdots\\
	&&&&&&\phantom{P}&&P_3\arrow[dl]& P_1\arrow[l] & P_2 \arrow[l] & S_{34}\arrow[l,left hook->]\\
	&&&&P_3\arrow[dl]&\arrow[l]\cdots&\cdots& P_2\arrow[l]&P_4\arrow[l]& P_1\arrow[l] & P_2 \arrow[l] & S_{34}\arrow[l,left hook->]\\
	S_4 & P_4\arrow[l,two heads]& P_1 \arrow[l]& P_2\arrow[l]&&&&&&&&P_3\arrow[dl]&P_1\arrow[l]&P_2\arrow[l]& S_{34}\arrow[l,left hook->]\\
	&&&&P_4\arrow[ul]&\arrow[l]\cdots&\cdots&P_2\arrow[l]&P_3\arrow[l]& P_1\arrow[l] & P_2 \arrow[l] & P_4\arrow[l]& P_1\arrow[l] & P_2 \arrow[l] & S_{34}\arrow[l,left hook->]\\
	&&&&&&\phantom{P}&&P_4\arrow[ul]& P_1\arrow[l] & P_2 \arrow[l] & P_3\arrow[l]& P_1\arrow[l] & P_2 \arrow[l] & S_{34}\arrow[l,left hook->]\\
	&&&&&&&&&&&P_4\arrow[ul]& P_1\arrow[l] & P_2 \arrow[l] & S_{34}\arrow[l,left hook->]\\
	&&&&&&\vdots&\vdots&\vdots&\vdots&\vdots&\vdots&\vdots&\vdots&\vdots\\
	&&&&&&&&&&&P_3\arrow[dl]& P_1\arrow[l] & P_2 \arrow[l] & S_{34}\arrow[l,left hook->]\\
	&&&&&&\cdots&P_2\arrow[l]&P_3\arrow[l]& P_1\arrow[l] & P_2 \arrow[l] & P_4\arrow[l]& P_1\arrow[l] & P_2 \arrow[l] & S_{34}\arrow[l,left hook->]\\
	&&&&&&&&P_4\arrow[ul]& P_1\arrow[l] & P_2 \arrow[l] & P_3\arrow[l]& P_1\arrow[l] & P_2 \arrow[l] & S_{34}\arrow[l,left hook->]\\
	&&&&&&&&&\phantom{P}&&P_4\arrow[ul]& P_1\arrow[l] & P_2 \arrow[l] & S_{34}\arrow[l,left hook->]\\
	\text{\footnotesize$-1$}&\text{\footnotesize$0$}&\text{\footnotesize$1$}&\text{\footnotesize$2$}&\text{\footnotesize$3$}&\text{\footnotesize$\cdots$}&\text{\footnotesize$\cdots$}&\text{\footnotesize$3\ell-1$}&\text{\footnotesize$3\ell$}&\text{\footnotesize$3\ell+1$}&\text{\footnotesize$\cdots$}&\text{\footnotesize$3k$}&\text{\footnotesize$3k+1$}&\text{\footnotesize$3k+2$}&\text{\footnotesize$3k+3$}
\end{tikzcd}
\end{center}
For clarity, we have indicated the degree of each term below the diagram. The first $3\ell$ columns coincide with the minimal projective resolution of $S_4$. This means that for $\ell \leq j < k$, the complex $X^\ell$ has $2^\ell$ indecomposable direct summands in degrees $3j+1$ and $3j+2$. Likewise, $X^\ell$ has $2^{\ell+1}$ indecomposable direct summands in degree $3j+3$ (with each term $S_{34}$ counted twice).

As in the example, Proposition \ref{prop:syzygy} implies that $Z^3_\ell$ is a direct summand of $\Omega^3_{C^b(A)} X^\ell$. Indeed, this is precisely the direct summand of $\Omega^{3\ell}X_k$ described in Corollary \ref{cor:formulaResults}(3).

Now, if $\ell < k$, then Lemma \ref{lem:smallSyzygy} implies that $\Omega^{3(k-\ell-1)}_{C^b(A)} Z_\ell^3 \cong (Z^3_\ell \oplus Z^4_\ell)^{2^{k-\ell-2}}$ is a direct summand of $\Omega^{3k}_{C^b(A)} X_k$. If $\ell = k-1$, then we have $\Omega^{3}_{C^b(A)}X^{k-1} \cong Z_{k-1}^{3} \oplus X^k \cong Z_{k-1}^{3} \oplus Z_{k}^{4}$, which are the final direct summands of $\Omega^{3k}_{C^b(A)} X_k$.
\end{proof}

\noindent We are now ready to prove Proposition \ref{prop:phidiminfinite}.

\begin{proof}[Proof of Proposition \ref{prop:phidiminfinite}]\

	Let $k \geq 1$. Recall from Lemma \ref{lem:commutes} that syzygy and wrapping commute. Thus by Lemma \ref{lem:bigSyzygy}, we have:
	\begin{eqnarray*}
		\Omega_{C^3(A)}^{3k+1}WX_k &\cong& WZ_k^1 \oplus WZ_{k-1}^1 \oplus \left(\bigoplus_{j = 0}^{k-2} (WZ_j^1\oplus WZ_j^1)^{2^{k-j-2}}\right) \cong \Omega_{C^3(A)}^{3k+1}WY_k\\
		\Omega_{C^3(A)}^{3k}WX_k &\cong& WZ_k^4 \oplus WZ_{k-1}^3 \oplus \left(\bigoplus_{j = 0}^{k-2} (WZ_j^3\oplus WZ_j^4)^{2^{k-j-2}}\right)\\
		\Omega_{C^3(A)}^{3k}WY_k &\cong& WZ_k^3 \oplus WZ_{k-1}^4 \oplus \left(\bigoplus_{j = 0}^{k-2} (WZ_j^3\oplus WZ_j^4)^{2^{k-j-2}}\right).
	\end{eqnarray*}
	This means $\Omega_{C^3(A)}^{3k}WX_k \cong \Omega_{C^3(A)}^{3k}WY_k$ if and only if $WZ_k^3 \oplus WZ_{k-1}^4 \cong WZ_k^4 \oplus WZ_{k-1}^3$. However, $WZ_k^3 \ncong WZ_k^4$ since $S_3$ is a direct summand of $(WZ_k^3)_{[-1]}$ but not $(WZ_k^4)_{[-1]}$. Likewise, $WZ_k^3 \ncong WZ_{k-1}^3$ since $(S_3\oplus S_4)^{2^k}$ is a direct summand of $(WZ_k^3)_{[-1]}$ but not $(WZ_{k-1}^3)_{[-1]}$. Since each of these 3-cyclic chain complexes is indecomposable by Lemma \ref{lem:indecomposable2}, we conclude that $WZ_k^3 \oplus WZ_{k-1}^4 \ncong WZ_k^4 \oplus WZ_{k-1}^3$.
\end{proof}

\section{Conclusions and Directions for Future Research}
\noindent As stated in the introduction, $\mathsf{fin.dim}\left(A\otimes_K A_3^{CT}\right) = 0$ as a consequence of \cite[Thorem 16]{eilenberg_dimension}. In particular, this means that for $X_k$ and $Y_k$ defined in Section \ref{sec:computation}, we must have $\mathsf{pd}_{C^3(A)} WX_k = \infty = \mathsf{pd}_{C^3(A)} WY_k$.

We can also see this directly, as both $X_k$ and $Y_k$ are truncated projective resolutions of modules of infinite projective dimension. Thus, for any integer $n$, the supports of $\Omega_{C^3(A)}^n WX_k$ and $\Omega^n_{C^3(A)}WY_k$ will contain non-projective $A$-modules in degree $[-1]$. Similar reasoning further shows that truncated projective resolutions could only be used to show a category of 3-periodic chain complexes has infinite finitistic dimension if the original algebra was already of infinite finitistic dimension. This is consistent with \cite[Thorem 16]{eilenberg_dimension}. It is nevertheless possible that a variant of this construction could give a counterexample to the finitistic dimension conjecture.

Our immediate goal moving forward is to better understand the class of algebras with infinite $\phi$-dimension. If there is a counterexample to the finitistic dimension conjecture, it will necessarily come from this class. Thus we wish to determine a more nonspecific method to determine whether an algebra is of infinite $\phi$-dimension. We would also like to find a less complicated example of such an algebra.

One possibility is the commutative algebra $K[x,y,z]/(x^2,y^2,z^2,xy)$. We believe this algebra is also of infinite $\phi$-dimension because is it isomorphic to the algebra $K[x,y]/(x^2,y^2,xy)\otimes_K K[z]/(z^2)$. Thus we can view the module category as a category of `1-periodic chain complexes'. However, our proof of Lemma \ref{lem:indecomposable2} does not generalize to this example, meaning a more subtle argument is needed. Moreover, this algebra cannot possibly be a counterexample to the finitistic dimension conjecture since it is commutative.

\section*{Acknowledgements}
Both authors are thankful to Gordana Todorov, who is working with them on the larger study of amalgamation, for numerous meaningful conversations and support. They would also like to thank Kaveh Mousavand, who asked them whether they could generalize the result of \cite{hanson_resolution} to arbitrary monomial relation algebras, starting the investigation that led to this paper. The second author is thankful to Daniel \'Alvarez-Gavela for their collaboration on the use of amalgamation to describe examples and invariants in contact topology (see \cite{alvarez_turaev}), which motivated the larger study of amalgamation. Moreover, the second author thanks An Huang for sharing the reference \cite{arkani_scattering} with him, which laid the ground work for the connection between \cite{alvarez_turaev}, \cite{hanson_resolution}, and this paper. Lastly, the authors are thankful to Rene Marczinzik for pointing out that $\mathsf{fin.dim}\left(A\otimes_K A_3^{CT}\right) = 0$ as a consequence of \cite[Thorem 16]{eilenberg_dimension} and to both Liang Chen and an anonymous referee for pointing out a errors in Lemma \ref{lem:bigSyzygy}.


\bibliography{Hanson_Igusa_MathZ_May2021.bib}
\end{document}